\documentclass[12pt]{article}
\usepackage{amsmath}
\usepackage{amssymb} 
\begin{document}

\newtheorem{df}{Definition}
\newtheorem{thm}{Theorem}
\newtheorem{lm}{Lemma}
\newtheorem{pr}{Proposition}
\newtheorem{co}{Corollary}
\newtheorem{re}{Remark}
\newtheorem{note}{Note}
\newtheorem{claim}{Claim}

\begin{center}
\begin{Large}\begin{bf}Characterization of quasi-Banach spaces which coarsely embed into a Hilbert space\end{bf}\end{Large}

\indent

N.~Lovasoa~RANDRIANARIVONY \footnote{Supported in part by 
NSF 0200690 and Texas Advanced Research Program 010366-0033-20013.

This paper represents a portion of the author's dissertation being prepared at 
Texas A\&M University under the direction of W.~B.~Johnson.}

\end{center}

\indent

A (not necessarily continuous) map $f$ between two metric spaces $(X,d)$ and $(Y, \delta )$ is called a 
{\sl coarse embedding\/} (see \cite[7.G]{G}) if there exist two non-decreasing functions $\varphi _1:[0,\infty)
\rightarrow [0,\infty )$ and  $\varphi _2:[0,\infty) \rightarrow [0,\infty )$ such that 
 
\begin{enumerate}
\item $\varphi_1(d(x,y))\leq \delta (f(x),f(y)) \leq \varphi_2(d(x,y))$
\item $\varphi _1(t) \rightarrow \infty$ as $t\rightarrow \infty$.
\end{enumerate}

It was proved in \cite{JR}  that $\ell_p$ does not coarsely embed into a Hilbert space when $p>2$.  The present article is a
strengthening of that result by giving a full characterization of quasi-Banach spaces that coarsely embed into a
Hilbert space.  This result, as well as its proof, mirrors the theorem in \cite{AMM} which characterizes spaces
that uniformly embed into a Hilbert space.
 The combination of Theorem 1 and the theorem in 
\cite{AMM} yield that a quasi-Banach space uniformly embeds into a Hilbert space if and only if it coarsely
embeds into a Hilbert space.  This is counterintuitive in that a uniform embedding gives information only on
small distances while a coarse embedding gives information only on large distances.

\indent

\begin{thm} A quasi-Banach space $X$ coarsely embeds into a Hilbert space if and only  if there is a probability
space $(\Omega, \mathcal{B}, \mu)$ such that $X$ is linearly isomorphic to a subspace of $L_0(\mu)$.\end{thm}

\indent

To prove Theorem 1, we use a result essentially contained in \cite{JR}  as is recalled in the following Proposition: 

\indent

\begin{pr} Let $X$ be a quasi-Banach space which coarsely embeds into a Hilbert space.  Then there exists on $X$ a
continuous negative definite function $g$ which satisfies $g(0)=0$ and $\phi_1(\|x\|)\leq g(x) \leq
\|x\|^{2\alpha}$ where $\phi_1 :[0,\infty) \rightarrow [0,\infty)$ is a nondecreasing function satisfying
$\phi_1(t) \rightarrow \infty$ as $t \rightarrow \infty$, and $\alpha >0$.   \end{pr}

\indent

\begin{it}Proof:\end{it} Steps 0, 1, 2 (and a piece of Step 3 on the continuity of $g$) in \cite{JR} extend to the case of quasi-Banach spaces.

\indent

\begin{it}Proof of Theorem 1:\end{it}

Let $X$ be a quasi-Banach space.  A theorem of Aoki-Rolewicz gives an equivalent quasi-norm $\|\cdot \|$ on $X$ which is also $p$-additive for some $0<p\leq 1$, i.e. $\|x+y\|^p \leq \|x\|^p + \|y\|^p$ for all $x, y$ in $X$.  In particular $X$ under this norm has type $p$.

Say $X$ is linearly isomorphic to a subspace of $L_0(\mu)$ for some probability space  $(\Omega, \mathcal{B}, \mu )$.  Then since $X$ has type $p$, $X$ is isomorphic to a subspace of $L_r(\mu)$ for every $r<p$.  (See \cite {BL}, Theorem 8.15.)

Now since $r<2$, Nowak's result \cite{N} implies that $X$ coarsely embeds into a Hilbert space.

\indent

Conversely, let $X$ be a quasi-Banach space which coarsely embeds into a Hilbert space.  Let $g$ be the negative
definite function on $X$ given by Proposition 1, and let $f$ be the continuous positive definite function given
by
$f=e^{-g}$.  Use Proposition 8.7 in \cite{BL} to get a probability space $(\Omega, \mathcal{B}, \mu)$ and a
continuous linear operator $U: X\rightarrow L_0(\mu)$ such that the characteristic function $\mathbb{E}\exp
(itUx)$ of $Ux$ is equal to $f(tx)$ for every $x \in X$ and $t \in \mathbb{R}$.  We show that $U$ is an
isomorphism into.

Let $(x_n)_n$ be a sequence in $X$ such that $U(x_n)\rightarrow 0$ in $L_0(\mu)$, i.e. in measure.  Then
$f(tx_n)=\mathbb{E}(\exp (itUx_n)) \rightarrow 1$ for each fixed $t$ in $\mathbb{R}$.  

If $(x_n)_n$ does not converge to $0$, then by passing to a subsequence we can assume without loss of generality that $\|x_n\| \geq
\epsilon$ for all $n$ and for some $\epsilon>0$.  But since $\phi_1$ is nondecreasing, we get for every $t>0$:

$$e^{-\phi_1(t\|x_n\|)} \leq e^{-\phi_1(t\epsilon)}.$$

\newpage

Since $\phi_1(s)\to \infty$ as $s\to \infty$, we can 
pick  $t_0>0$ so  that $e^{-\phi_1(t_0\epsilon)}<\frac{1}{2}$.  For that $t_0$, we have for every
$n$:

$$f(t_0x_n)\leq e^{-\phi_1(t_0\epsilon)} < \frac{1}{2}.$$

This contradicts  the fact that $f(t_0x_n) \rightarrow 1$.

Thus  $x_n \rightarrow 0$, and hence $U$ is one-to-one and its inverse is continuous.

\indent

\noindent \begin{bf}Acknowledgement:\end{bf}  I would like to thank W.~B.~Johnson for the constructive conversations that led to the production of this article.

\bigskip

\bigskip

\noindent    N.~Lovasoa Randrianarivony\newline
             Department of Mathematics\newline
             Texas A\&M University\newline
             College Station, TX, USA\newline
             E-mail: nirina@math.tamu.edu

\end{document}